\documentclass[12pt]{article}
\usepackage{amssymb,amsmath}
\usepackage{cases}
\usepackage{amsfonts}
\usepackage{color,xcolor}
\usepackage[left=2.0cm,right=2.0cm,top=2.0cm,bottom=2.0cm]{geometry}
\usepackage[colorlinks,citecolor=blue,urlcolor=blue]{hyperref}
\usepackage[pagewise]{lineno}

\newtheorem{theorem}{Theorem}[section]

\newtheorem{lemma}{Lemma}[section]
\newtheorem{proposition}{Proposition}[section]

\newtheorem{definition}{Definition}[section]
\newtheorem{remark}{Remark}[section]

\newcommand{\bal}{\begin{align}}
\newcommand{\bbal}{\begin{align*}}
\newcommand{\beq}{\begin{equation}}
\newcommand{\eeq}{\end{equation}}
\newcommand{\bca}{\begin{cases}}
\newcommand{\eca}{\end{cases}}

\newcommand{\pa}{\partial}

\newcommand{\dd}{\mathrm{d}}

\newcommand{\R}{\mathbb{R}}

\newcommand{\Z}{\mathbb{Z}}

\newcommand{\bi}{\Big}

\begin{document}
\title{Non-uniform continuity of the Fokas-Olver-Rosenau-Qiao
equation in Besov spaces}

\author{Xing Wu\footnote{E-mail:ny2008wx@163.com}\\
\small \it College of Information and Management Science,
Henan Agricultural University,\\
\small Zhengzhou, Henan, 450002, China}

\date{}

\maketitle\noindent{\hrulefill}

{\bf Abstract:} In this paper, we prove that the solution map of Fokas-Olver-Rosenau-Qiao
equation (FORQ) is not uniformly continuous
on the initial data in Besov spaces. Our result  extends the previous non-uniform continuity in Sobolev spaces (Nonlinear Anal., 2014) \cite{Himonas 2014} to Besov spaces and is consistent with  the present work (J. Math. Fluid Mech., 2020) \cite{1Li 2020} on Novikov equation up to some coefficients when dropping the extra term $(\partial_xu)^3$ in FORQ.


\vskip0mm\noindent{\hrulefill}

\section{Introduction}

In this paper, we are concerned with the following Fokas-Olver-Rosenau-Qiao equation (FORQ)
\begin{eqnarray}\label{eq1}
        \left\{\begin{array}{ll}
         u_t-u_{xxt}+3u^2u_x-u_x^3-4uu_xu_{xx}+2u_xu_{xx}^2-u^2u_{xxx}+u_x^2u_{xxx}=0, ~~t>0, ~x\in \mathbb{R},\\
          u(0, x)=u_0, ~~x\in \mathbb{R}.\end{array}\right.
        \end{eqnarray}
Eq. (\ref{eq1}) written in a slightly different form was first derived
by Fokas \cite{Fokas 1995} as an integrable generalisation of the modified KdV equation. Soon after, Fuchssteiner \cite{Fuchssteiner 1996}
and Olver-Rosenau \cite{Olver 1996} independently obtained similar versions of this equation by performing a simple explicit algorithm based on the bi-Hamiltonian representation of the classically integrable system. Several years later, the concise form written above was recovered  by Qiao \cite{Qiao 2006} from the two-dimensional Euler equations by using an approximation procedure.

The entire integrable hierarchy related
to the FORQ equation was proposed by Qiao \cite{Qiao 2007}.  It also has bi-Hamiltonian structure, which was first derived in \cite{Olver 1996} and then in \cite{Qiao 2006},  admits Lax pair \cite{Qiao 2006} and peakon travelling wave solutions that  are orbitally stable \cite{Gui 2013, Qu 2013, Liu 2014}. For more discussion about Lax integrability and peakon solutions of FORQ we refer to \cite{Chang 2018}, where this equation is also referred as the modified Camassa-Holm equation.

The local well-posedness and
ill-posedness for the Cauchy problem of the FORQ equation (\ref{eq1}) in Sobolev spaces and Besov spaces were studied in the series of papers  \cite{Himonas 2014, 1Himonas 2014, Himonas 2019, Fu 2013}. It was showed by  Himonas-Mantzavinos \cite{Himonas 2014} that the FORQ is well-posed in Sobolev space $H^s$ with $s>\frac{5}{2}$ in the sense of  Hadamard. Fu et al.\cite{Fu 2013} established the local well-posedness in Besov space  $B_{p, r}^s$ with  $s>\max\{2+\frac{1}{p}, \frac{5}{2}\}$, $1\leq p, r\leq \infty$.
 After the non-uniform dependence  for some dispersive equations was studied by Kenig et al. \cite{Kenig 2001}, the issue of non-uniform continuity of solutions on initial data has attracted much more attention, such as on classical Camassa-Holm equation \cite{Himonas 2005, Himonas 2009, Himonas 2010, Li 2020} and on famous Novikov equation \cite{Himonas 2012, 1Li 2020}. It was further proved in \cite{Himonas 2014} that the dependence on initial data is sharp, i.e. the data-to-solution map is continuous but not uniformly continuous.

For studying the non-uniform continuity of the FORQ equation, it is
more convenient to express (\ref{eq1}) in the following equivalent nonlocal form
\begin{equation}\label{eq2}
\begin{cases}
u_t+u^2\pa_xu=\frac{1}{3}(\partial_xu)^3-\frac{1}{3}(1-\partial_x^2)^{-1}[(\partial_xu)^3]-
\partial_x(1-\partial_x^2)^{-1}[\frac{2}{3}u^3+u(\partial_xu)^2],\\
u(0,x)=u_0,\;\qquad  t>0, ~x\in \mathbb{R}.
\end{cases}
\end{equation}
When removing $(\partial_xu)^3$ from (\ref{eq2}), (\ref{eq2}) becomes the following Novikov equation up to some coefficients
\begin{equation}\label{eq3}
\begin{cases}
u_t+u^2\pa_xu=-\frac{1}{3}(1-\partial_x^2)^{-1}[(\partial_xu)^3]-
\partial_x(1-\partial_x^2)^{-1}[\frac{2}{3}u^3+u(\partial_xu)^2],\\
u(0,x)=u_0,\;\qquad  t>0, ~x\in \mathbb{R}.
\end{cases}
\end{equation}
Recently, Li, Yu and Zhu \cite{1Li 2020} have proved that the solution map of Novikov equation is not uniformly continuous dependence on the initial data in the Besov spaces $B_{p, r}^s(\mathbb{R})$, $s>\max\{1+\frac{1}{p}, \frac{3}{2}\}$, $1\leq p, r\leq \infty$.  It is noticed that well-posedness for Novikov equation holds for $s>\max\{1+\frac{1}{p}, \frac{3}{2}\}$ while well-posedness for FORQ holds
for $s>\max\{2+\frac{1}{p}, \frac{5}{2}\}$. This difference between the well-posedness index of these equations may be explained by the presence of the extra term $(\partial_xu)^3$ in FORQ, which is not quasi-linear and is absent from Novikov equation \cite{Himonas 2014}.

Up to the present, there is no result for the non-uniform continuous dependence of FORQ in Besov space and it seems more difficult due to the presence of the extra term $(\partial_xu)^3$, which elevates two regularities, compared with Novikov equation and the method developed for the Novikov equation in \cite{1Li 2020} will make it more complex. In this paper, we will follow a different route to bypass this problem. Firstly, consider a new system satisfied by $(1-\partial_x)u:\triangleq v$. Secondly, for any bounded set $v_0$ in working space, the corresponding solution $\mathbf{S}_t(v_0)$ can be approximated by a function of first degree of time $t$ with convective term and nonlocal term being coefficients. With suitable choice of initial data, the difference between two solutions will produce a term from  convective term which will not be  small for small time, and thus we obtain the non-uniform continuous dependence of FORQ.
These will be described in more detail later.

Now, we state our main result.
\begin{theorem}\label{the1.1} Let $s>\max\{2+\frac{1}{p}, \frac{5}{2}\}$, $1\leq p, r\leq \infty$. The solution map $u_0\rightarrow \mathbf{S}_t(u_0)$ of the initial value problem (\ref{eq2}) is not
uniformly continuous from any bounded subset of  $B_{p, r}^s(\mathbb{R})$ into $\mathcal{C}([0, T];  B_{p, r}^s(\mathbb{R}))$. More precisely,
there exist two sequences $u^{1, n}(0, x)$ and $u^{2, n}(0, x)$ such that
\begin{eqnarray*}
        \|u^{1, n}(0, x), u^{2, n}(0, x)\|_{B_{p, r}^s}\lesssim 1,\qquad  \qquad \qquad \lim_{n\rightarrow \infty} \|u^{1, n}(0, x)-u^{2, n}(0, x)\|_{B_{p, r}^s}=0,
        \end{eqnarray*}
but
\begin{eqnarray*}
      \liminf_{n\rightarrow \infty} \|\mathbf{S}_t(u^{1, n}(0, x))-\mathbf{S}_t(u^{2, n}(0, x))\|_{B_{p, r}^s}\gtrsim t, \qquad t\in [0, T_0],
        \end{eqnarray*}
with small positive time $T_0$ for $T_0\leq T$.
\end{theorem}

\begin{remark}\label{rem1}
Since $B_{2, 2}^s=H^s$ for any $s\in \mathbb{R}$, our result  extends the previous non-uniform continuity in Sobolev spaces \cite{Himonas 2014} to Besov spaces.
\end{remark}

\begin{remark}\label{rem2}
When dropping the extra term $(\partial_xu)^3$ in FORQ, we can get the same result on Novikov equation up to some coefficients, which is consistent with  the present work \cite{1Li 2020} on Novikov equation. The method we use in proving Theorem \ref{the1.1} is different from \cite{1Li 2020}  and is more general.
\end{remark}

{\bf Notations}:  Given a Banach space $X$, we denote the norm of a function on $X$ by $\|\|_{X}$, and \begin{eqnarray*}
\|\cdot\|_{L_T^\infty(X)}=\sup_{0\leq t\leq T}\|\cdot\|_{X}.
\end{eqnarray*} The symbol
$A\lesssim B$ means that there is a uniform positive constant $C$ independent of $A$ and $B$ such that $A\leq CB$.

\section{Littlewood-Paley analysis}

In this section, we will review the definition of Littlewood-Paley decomposition and nonhomogeneous Besov space, and then list some useful properties. For more details, the readers can refer to \cite{Bahouri 2011}.

There exists a couple of smooth functions $(\chi,\varphi)$ valued in $[0,1]$, such that $\chi$ is supported in the ball $\mathcal{B}\triangleq \{\xi\in\mathbb{R}^d:|\xi|\leq \frac 4 3\}$, $\varphi$ is supported in the ring $\mathcal{C}\triangleq \{\xi\in\mathbb{R}^d:\frac 3 4\leq|\xi|\leq \frac 8 3\}$ and $\varphi\equiv 1$ for $\frac{4}{3}\leq |\xi| \leq \frac{3}{2}$. Moreover,
$$\forall\,\, \xi\in\mathbb{R}^d,\,\, \chi(\xi)+{\sum\limits_{j\geq0}\varphi(2^{-j}\xi)}=1,$$
$$\forall\,\, \xi\in\mathbb{R}^d\setminus\{0\},\,\, {\sum\limits_{j\in \mathbb{Z}}\varphi(2^{-j}\xi)}=1,$$
$$|j-j'|\geq 2\Rightarrow\textrm{Supp}\,\ \varphi(2^{-j}\cdot)\cap \textrm{Supp}\,\, \varphi(2^{-j'}\cdot)=\emptyset,$$
$$j\geq 1\Rightarrow\textrm{Supp}\,\, \chi(\cdot)\cap \textrm{Supp}\,\, \varphi(2^{-j}\cdot)=\emptyset.$$
Then, we can define the nonhomogeneous dyadic blocks $\Delta_j$ and nonhomogeneous low frequency cut-off operator $S_j$ as follows:
$$\Delta_j{u}= 0,\,\, if\,\, j\leq -2,\quad
\Delta_{-1}{u}= \chi(D)u=\mathcal{F}^{-1}(\chi \mathcal{F}u),$$
$$\Delta_j{u}= \varphi(2^{-j}D)u=\mathcal{F}^{-1}(\varphi(2^{-j}\cdot)\mathcal{F}u),\,\, if \,\, j\geq 0,$$
$$S_j{u}= {\sum\limits_{j'=-\infty}^{j-1}}\Delta_{j'}{u}.$$

\begin{definition}[\cite{Bahouri 2011}]\label{de2.3}
Let $s\in\mathbb{R}$ and $1\leq p,r\leq\infty$. The nonhomogeneous Besov space $B^s_{p,r}(\R^d)$ consists of all tempered distribution $u$ such that
\begin{align*}
||u||_{B^s_{p,r}(\R^d)}\triangleq \Big|\Big|(2^{js}||\Delta_j{u}||_{L^p(\R^d)})_{j\in \Z}\Big|\Big|_{\ell^r(\Z)}<\infty.
\end{align*}
\end{definition}

In the following, we list some basic lemmas and properties about Besov space which will be frequently used in proving our main result.

\begin{lemma}(\cite{Bahouri 2011})\label{lem2.1}
 (1) Algebraic properties: $\forall s>0,$ $B_{p, r}^s(\mathbb{R}^d)$ $\cap$ $L^\infty(\mathbb{R}^d)$ is a Banach algebra. $B_{p, r}^s(\mathbb{R}^d)$ is a Banach algebra $\Leftrightarrow B_{p, r}^s(\mathbb{R}^d)\hookrightarrow L^\infty(\mathbb{R}^d)\Leftrightarrow s>\frac{d}{p}$ or $s=\frac{d}{p},$ $r=1$.\\
 (2) For any $s>0$ and $1\leq p,r\leq\infty$, there exists a positive constant $C=C(d,s,p,r)$ such that
$$\|uv\|_{B^s_{p,r}(\mathbb{R}^d)}\leq C\Big(\|u\|_{L^{\infty}(\mathbb{R}^d)}\|v\|_{B^s_{p,r}(\mathbb{R}^d)}+\|v\|_{L^{\infty}(\mathbb{R}^d)}\|u\|_{B^s_{p,r}(\mathbb{R}^d)}\Big).$$
(3) Let $m\in \mathbb{R}$ and $f$ be an $S^m-$ multiplier (i.e., $f: \mathbb{R}^d\rightarrow \mathbb{R}$ is smooth and satisfies that $\forall \alpha\in \mathbb{N}^d$, there exists a constant $\mathcal{C}_\alpha$ such that $|\partial^\alpha f(\xi)|\leq \mathcal{C}_\alpha(1+|\xi|)^{m-|\alpha|}$ for all $\xi \in \mathbb{R}^d$). Then the operator $f(D)$ is continuous from $B_{p, r}^s(\mathbb{R}^d)$ to $B_{p, r}^{s-m}(\mathbb{R}^d)$.\\
(4) For any $s\in \mathbb{R}$, $(1-\partial_x)^{-1}$ is an isomorphic mapping from $B_{p, r}^{s-1}(\mathbb{R}^d)$ into $B_{p, r}^s(\mathbb{R}^d)$.
\end{lemma}

\begin{lemma}(\cite{Bahouri 2011, Li-Yin2})\label{lem2.2}
Let $1\leq p,r\leq \infty$ and
$
\sigma> - \min\{\frac{1}{p}, 1-\frac{1}{p}\}.
$
There exists a constant $C=C(p,r,\sigma)$ such that for any smooth solution to the following linear transport equation:
\begin{equation*}
\quad \partial_t f+v\pa_xf=g,\quad \; f|_{t=0} =f_0.
\end{equation*}
We have
\begin{align}\label{ES2}
\sup_{s\in [0,t]}\|f(s)\|_{B^{\sigma}_{p,r}(\mathbb{R})}\leq Ce^{CV_{p}(v,t)}\Big(\|f_0\|_{B^\sigma_{p,r}(\mathbb{R})}
+\int^t_0\|g(\tau)\|_{B^{\sigma}_{p,r}(\mathbb{R})}\dd \tau\Big),
\end{align}
with
\begin{align*}
V_{p}(v,t)=
\begin{cases}
\int_0^t \|\nabla v(s)\|_{B^{\frac{1}{p}}_{p,\infty}(\mathbb{R})\cap L^\infty(\mathbb{R})}\dd s,&\ \ \mathrm{if} \; \sigma<1+\frac{1}{p},\\
\int_0^t \|\nabla v(s)\|_{B^{\sigma}_{p,r}(\mathbb{R})}\dd s,&\ \ \mathrm{if} \; \sigma=1+\frac{1}{p} \mbox{ and } r>1,\\
\int_0^t \|\nabla v(s)\|_{B^{\sigma-1}_{p,r}(\mathbb{R})}\dd s, &\ \ \mathrm{if} \;\sigma>1+\frac{1}{p}\ \mathrm{or}\ \{\sigma=1+\frac{1}{p} \mbox{ and } r=1\}.
\end{cases}
\end{align*}
\end{lemma}

\section{Reformulation of the System}
\setcounter{equation}{0}
Due to the presence of the extra term $(\partial_xu)^3$ in FORQ,  it seems difficult to deal with Eq. (\ref{eq2}) directly. Therefore, we shall first differentiate FORQ with respect to $x$ and then simplify the resulting expression, we obtain
\begin{eqnarray}\label{eq3.1}
         \partial_t(\partial_xu)&=&(\partial_xu)^2\partial_x^2u-2u(\partial_xu)^2-u^2\partial_x^2u+[\frac{2}{3}u^3
         +u(\partial_xu)^2]\nonumber \\
        &\;&-(1-\partial_x^2)^{-1}\partial_x[\frac{1}{3}(\partial_xu)^3]-(1-\partial_x^2)^{-1}[\frac{2}{3}u^3
         +u(\partial_xu)^2].
        \end{eqnarray}
Let $v=(1-\partial_x)u$, we have from (\ref{eq2}) and (\ref{eq3.1}) that
\begin{eqnarray}\label{eq3.2}
        \left\{\begin{array}{ll}
        \partial_tv=(v^2-2uv)\partial_xv-\frac{1}{3}u^3-\frac{1}{3}v^3-\Phi_1(v)-\Phi_2(v),\\
          u=(1-\partial_x)^{-1}v,\\
          v(0, x)=(1-\partial_x)u_0(x)\triangleq v_0, \end{array}\right.
        \end{eqnarray}
where the nonlocal terms $\Phi_1(v), \Phi_2(v)$ are defined by
\begin{eqnarray*}
        \Phi_1(v)=(1-\partial_x^2)^{-1}[\frac{8}{3}u^3-\frac{1}{3}v^3-3u^2v], \;\;\; \Phi_2(v)=\partial_x(1-\partial_x^2)^2[\frac{1}{3}v^3-u^2v].
        \end{eqnarray*}
Since $(1-\partial_x)^{-1}$ is an isomorphic mapping from $B_{p, r}^{s-1}(\mathbb{R})$ into $B_{p, r}^s(\mathbb{R})$, the non-uniform continuous dependence of $u$ in $B_{p, r}^s$ then can be transformed into that of $v$ in $B_{p, r}^{s-1}.$

\section{Non-uniform continuous dependence}
\setcounter{equation}{0}

In this section, we will give the proof of Theorem \ref{the1.1}. However, as explained above, we will directly consider Eq. (\ref{eq3.2}) satisfied by $v$.

Firstly, we establish the estimates of the difference between the solution $\mathbf{S}_t(v_0)$ and  initial data $v_0$ in different Besov norms. That is

\begin{proposition}\label{pro2}
Assume that $||v_0||_{B^{s-1}_{p,r}}\lesssim 1$. Under the assumptions of Theorem \ref{the1.1}, we have
\bbal
&||\mathbf{S}_{t}(v_0)-v_0||_{B^{s-2}_{p,r}}\lesssim t||v_0||^2_{B^{s-2}_{p,r}}||v_0||_{B^{s-1}_{p,r}},
\\&||\mathbf{S}_{t}(v_0)-v_0||_{B^{s-1}_{p,r}}\lesssim t\big(||v_0||^3_{B^{s-1}_{p,r}}+||v_0||^2_{B^{s-2}_{p,r}}||v_0||_{B^{s}_{p,r}}\big),
\\&||\mathbf{S}_{t}(v_0)-v_0||_{B^{s}_{p,r}}\lesssim t\big(||v_0||^2_{B^{s-1}_{p,r}}||v_0||_{B^{s}_{p,r}}+||v_0||^2_{B^{s-2}_{p,r}}||v_0||_{B^{s+1}_{p,r}}\big).
\end{align*}
\end{proposition}
{\bf Proof}\quad For simplicity, denote $v(t)=\mathbf{S}_t(v_0)$. Firstly, according to the local well-posedness result \cite{Fu 2013, Himonas 2014}, there exists a positive time $T=T(||u_0||_{B^s_{p,r}},s,p,r)$ such that the solution $u(t)$ belongs to $\mathcal{C}([0, T];  B_{p, r}^s)$. Moreover, by Lemmas \ref{lem2.1}-\ref{lem2.2}, for all $t\in[0,T]$ and $\gamma\geq s-2$,  there holds
\bal\label{u-estimate}
||u(t)||_{B^\gamma_{p,r}}\leq C||u_0||_{B^\gamma_{p,r}}\;\;\text{or}\;\; ||v(t)||_{B^{\gamma-1}_{p,r}}\leq C||v_0||_{B^{\gamma-1}_{p,r}}.
\end{align}

Now we shall estimate the different Besov norms of the term $v(t)-v_0$, which can be bounded by $t$ multiplying the corresponding Besov norms of initial data $v_0$.

It follows by differential mean value theorem and the Minkowski inequality that
\bbal
||v(t)-v_0||_{B^{s-1}_{p,r}}
&\lesssim \int^t_0||\pa_\tau v||_{B^{s-1}_{p,r}} \dd\tau
\\&\lesssim \int^t_0||(v^2-2uv)\partial_xv||_{B^{s-1}_{p,r}} \dd\tau+ \int^t_0||\frac{1}{3}v^3||_{B^{s-1}_{p,r}} \dd\tau
\\&\quad+\int^t_0||\frac{1}{3}u^3||_{B^{s-1}_{p,r}} \dd\tau+\int^t_0||\Phi_1(v)||_{B^{s-1}_{p,r}} \dd\tau+\int^t_0||\Phi_2(v)||_{B^{s-1}_{p,r}} \dd\tau.
\end{align*}
Here, we shall only have to estimate $||(v^2-2uv)\partial_xv||_{B^{s-1}_{p,r}}$ and $||\frac{1}{3}v^3||_{B^{s-1}_{p,r}}$, since the other terms can be processed in a similar more relaxed way and have the same bound as $||v^3||_{B^{s-1}_{p,r}}$.

Using the fact that $B_{p, r}^{s-2}$ is an Banach algebra with $s-2>\max\{\frac{1}{p}, \frac{1}{2}\}$, together with the product estimates (2) in Lemma \ref{lem2.1}, one has
\begin{eqnarray*}
||v^3||_{B^{s-1}_{p,r}}&\lesssim& ||v||^3_{B^{s-1}_{p,r}},\\
||(v^2-2uv)\partial_xv||_{B^{s-1}_{p,r}}&\lesssim& ||(v^2-2uv)||_{L^\infty}||\partial_xv||_{B^{s-1}_{p,r}}+||(v^2-2uv)||_{B^{s-1}_{p,r}}||\partial_xv||_{L^\infty}\\
&\lesssim& ||v||^2_{B^{s-2}_{p,r}}||v||_{B^{s}_{p,r}}+||v||^3_{B^{s-1}_{p,r}},
\end{eqnarray*}
where we have used the relation $u=(1-\partial_x)^{-1}v$, and $(1-\partial_x)^{-1}$ is a  $S^{-1}-$ multiplier, which is continuous from $B_{p, r}^{s'-1}(\mathbb{R})$ to $B_{p, r}^{s'}(\mathbb{R})$, thus $||u||_{B^{s'}_{p,r}}\lesssim ||v||_{B^{s'-1}_{p,r}}$ for any $s'\in \mathbb{R}$.

Therefore, in view of (\ref{u-estimate}), for $t\in[0,T]$, we have
\bbal
||v(t)-v_0||_{B^{s-1}_{p,r}}
&\lesssim t||v||^2_{L_t^\infty(B^{s-2}_{p,r})}||v||_{L_t^\infty(B^{s}_{p,r})}+||v||^3_{L_t^\infty(B^{s-1}_{p,r})}\\
&\lesssim t(||v_0||^2_{B^{s-2}_{p,r}}||v_0||_{B^{s}_{p,r}}+||v_0||^3_{B^{s-1}_{p,r}}).
\end{align*}

Following the same procedure of estimate as above, we have
\begin{eqnarray*}
||v(t)-v_0||_{B^{s-2}_{p,r}}
&\lesssim& \int^t_0||\pa_\tau v||_{B^{s-2}_{p,r}} \dd\tau
\\&\lesssim& \int^t_0||(v^2-2uv)\partial_xv||_{B^{s-2}_{p,r}} \dd\tau+ \int^t_0||\frac{1}{3}v^3||_{B^{s-2}_{p,r}} \dd\tau
\\&\quad&+\int^t_0||\frac{1}{3}u^3||_{B^{s-2}_{p,r}} \dd\tau+\int^t_0||\Phi_1(v)||_{B^{s-2}_{p,r}} \dd\tau+\int^t_0||\Phi_2(v)||_{B^{s-2}_{p,r}} \dd\tau\\
&\lesssim & t||v||^2_{L_t^\infty(B^{s-2}_{p,r})}||v||_{L_t^\infty(B^{s-1}_{p,r})}\\
&\lesssim &t||v_0||^2_{B^{s-2}_{p,r}}||v_0||_{B^{s-1}_{p,r}},
\end{eqnarray*}
and
\begin{eqnarray*}
||v^3||_{B^{s}_{p,r}}&\lesssim& ||v^2||_{L^\infty}||v||_{B^{s}_{p,r}}+||v^2||_{B^{s}_{p,r}}||v||_{L^\infty}\\
&\lesssim& ||v||^2_{B^{s-1}_{p,r}}||v||_{B^{s}_{p,r}},\\
||(v^2-2uv)\partial_xv||_{B^s_{p,r}}&\lesssim& ||(v^2-2uv)||_{L^\infty}||\partial_xv||_{B^s_{p,r}}+||(v^2-2uv)||_{B^s_{p,r}}||\partial_xv||_{L^\infty}\\
&\lesssim& ||v||^2_{B^{s-2}_{p,r}}||v||_{B^{s+1}_{p,r}}+||v||^2_{B^{s-1}_{p,r}}||v||_{B^{s}_{p,r}},
\end{eqnarray*}
hence,
\begin{eqnarray*}
||v(t)-v_0||_{B^{s}_{p,r}}
&\lesssim& \int^t_0||\pa_\tau v||_{B^{s}_{p,r}} \dd\tau
\\&\lesssim& \int^t_0||(v^2-2uv)\partial_xv||_{B^{s}_{p,r}} \dd\tau+ \int^t_0||\frac{1}{3}v^3||_{B^{s}_{p,r}} \dd\tau
\\&\quad&+\int^t_0||\frac{1}{3}u^3||_{B^{s}_{p,r}} \dd\tau+\int^t_0||\Phi_1(v)||_{B^{s}_{p,r}} \dd\tau+\int^t_0||\Phi_2(v)||_{B^{s}_{p,r}} \dd\tau\\
&\lesssim & t(||v||^2_{L^\infty_t(B^{s-2}_{p,r})}||v||_{L^\infty_t(B^{s+1}_{p,r})}+||v||^2_{L^\infty_t(B^{s-1}_{p,r})}
||v||_{L^\infty_t(B^{s}_{p,r})})\\
&\lesssim &t(||v_0||^2_{B^{s-2}_{p,r}}||v_0||_{B^{s+1}_{p,r}}+||v_0||^2_{B^{s-1}_{p,r}}||v_0||_{B^{s}_{p,r}}),
\end{eqnarray*}
Thus, we finish the proof of Proposition \ref{pro2}.

With the different Besov norms  estimates  of $v-v_0$ at hand, we have the following core estimates, which implies that for any bounded initial data $v_0$ in $B^{s-1}_{p,r}$, the corresponding solution $S_t(v_0)$ can be approximated by $v_0+t(v_0^2-2u_0v_0)\partial_xv_0+\frac{1}{3}tv_0^3+t[\frac{1}{3}u_0^3+\Phi_1(v_0)+\Phi_2(v_0)]$ near $t=0$.
\begin{proposition}\label{pro3}
Assume that $||v_0||_{B^{s-1}_{p,r}}\lesssim 1$. Then under the assumptions of Theorem \ref{the1.1}, there holds
\bbal
||\mathbf{S}_{t}(v_0)-v_0-t\mathbf{v}_0||_{B^{s-1}_{p,r}}\lesssim t^{2}\big(||v_0||^3_{B^{s-1}_{p,r}}+||v_0||^2_{B^{s-2}_{p,r}}||v_0||_{B^{s}_{p,r}}
+||v_0||^{4}_{B^{s-2}_{p,r}}||v_0||_{B^{s+1}_{p,r}}\big),
\end{align*}
where $\mathbf{v}_0=(v_0^2-2u_0v_0)\partial_xv_0+\frac{1}{3}v_0^3+\frac{1}{3}u_0^3+\Phi_1(v_0)+\Phi_2(v_0).$
\end{proposition}
{\bf Proof}\quad Using differential mean value theorem and the Minkowski inequality, we first arrive at
\begin{eqnarray}\label{eq4.2}
||v(t)-v_0-t\mathbf{v}_0||_{B^{s-1}_{p,r}}
&\lesssim & \int^t_0||\pa_\tau v-\mathbf{v}_0||_{B^{s-1}_{p,r}} \dd\tau \nonumber\\
&\lesssim & \int^t_0||v^2\partial_xv-v_0^2\partial_xv_0||_{B^{s-1}_{p,r}} \dd\tau+ \int^t_0||2uv\partial_xv-2u_0v_0\partial_xv_0||_{B^{s-1}_{p,r}} \dd\tau\nonumber\\
 &\quad&+ \int^t_0||\frac{1}{3}v^3-\frac{1}{3}v_0^3||_{B^{s-1}_{p,r}} \dd\tau
+\int^t_0||\frac{1}{3}u^3-\frac{1}{3}u_0^3||_{B^{s-1}_{p,r}}\dd\tau \nonumber\\ &\quad&+\int^t_0||\Phi_1(v)-\Phi_1(v_0)||_{B^{s-1}_{p,r}} \dd\tau+\int^t_0||\Phi_2(v)-\Phi_2(v_0)||_{B^{s-1}_{p,r}} \dd\tau.
\end{eqnarray}
Using (3) in Lemma \ref{lem2.1}, it is sufficient to estimate $||v^2\partial_xv-v_0^2\partial_xv_0||_{B^{s-1}_{p,r}}$,  $||2uv\partial_xv-2u_0v_0\partial_xv_0||_{B^{s-1}_{p,r}}$ and $||v^3-v_0^3||_{B^{s-1}_{p,r}}$, since the other terms can be processed in a similar more relaxed way and have the same bound as $||v^3-v_0^3||_{B^{s-1}_{p,r}}$.

It should be noticed that according to (\ref{u-estimate}), $||v||_{B^{s-1}_{p,r}}\lesssim ||v_0||_{B^{s-1}_{p,r}}\lesssim 1$, which will be frequently used later.

Due to the fact that $B_{p, r}^{s-2}$ is an Banach algebra with $s-2>\max\{\frac{1}{p}, \frac{1}{2}\}$, combining with the product estimates (2) in Lemma \ref{lem2.1}, we get
\begin{eqnarray*}
  ||v^2\partial_xv-v_0^2\partial_xv_0||_{B^{s-1}_{p,r}}
  &=&||(v^2-v_0^2)\partial_xv+v_0^2(\partial_xv-\partial_xv_0)||_{B^{s-1}_{p,r}}\\
  &\lesssim & ||(v^2-v_0^2)\partial_xv||_{B^{s-1}_{p,r}}+||v_0^2(\partial_xv-\partial_xv_0)||_{B^{s-1}_{p,r}}\\
 &\lesssim& ||v^2-v_0^2||_{L^\infty}||\partial_xv||_{B^{s-1}_{p,r}}+  ||v^2-v_0^2||_{B^{s-1}_{p,r}} ||\partial_xv||_{L^\infty}\\
 &\;&+ ||v_0^2||_{L^\infty}||\partial_xv-\partial_xv_0||_{B^{s-1}_{p,r}}+  ||v_0^2||_{B^{s-1}_{p,r}} ||\partial_xv-\partial_xv_0||_{L^\infty}\\
 &\lesssim &||v-v_0||_{B^{s-2}_{p,r}}||v||_{B^s_{p,r}}+||v-v_0||_{B^{s-1}_{p,r}}
 +||v-v_0||_{B^s_{p,r}}||v_0||^2_{B^{s-2}_{p,r}},\\
||2uv\partial_xv-2u_0v_0\partial_xv_0||_{B^{s-1}_{p,r}}
&=&||2(u-u_0)v\partial_xv+2u_0(v\partial_xv-v_0\partial_xv_0)||_{B^{s-1}_{p,r}}\\
&\lesssim& ||u-u_0||_{B^{s-1}_{p,r}}||v\partial_xv||_{B^{s-1}_{p,r}}
+||u_0||_{B^{s-1}_{p,r}}||v\partial_xv-v_0\partial_xv_0||_{B^{s-1}_{p,r}}\\
&\lesssim& ||v-v_0||_{B^{s-2}_{p,r}}||v||_{B^s_{p,r}}+||v-v_0||_{B^{s-1}_{p,r}}
 +||v-v_0||_{B^s_{p,r}}||v_0||^2_{B^{s-2}_{p,r}},
\end{eqnarray*}
and
\begin{eqnarray*}
\|v^3-v_0^3\|_{B^{s-1}_{p,r}}=\|(v-v_0)(v^2+vv_0+v_0^2)\|_{B^{s-1}_{p,r}}
\lesssim ||v-v_0||_{B^{s-1}_{p,r}}.
\end{eqnarray*}
Taking the above estimates into (\ref{eq4.2}), which together with Proposition \ref{pro2} yield
\begin{eqnarray*}
||v(t)-v_0-t\mathbf{v}_0||_{B^{s-1}_{p,r}}
&\lesssim & \int_0^t(||v-v_0||_{B^{s-1}_{p,r}}+||v-v_0||_{B^{s-2}_{p,r}}||v||_{B^s_{p,r}}d\tau
+||v-v_0||_{B^s_{p,r}}||v_0||^2_{B^{s-2}_{p,r}})d\tau\\
&\lesssim & t^{2}\big(||v_0||^3_{B^{s-1}_{p,r}}+||v_0||^2_{B^{s-2}_{p,r}}||v_0||_{B^{s}_{p,r}}
+||v_0||^{4}_{B^{s-2}_{p,r}}||v_0||_{B^{s+1}_{p,r}}\big).
\end{eqnarray*}
Thus, we complete the proof of Proposition \ref{pro3}. $\Box$

Now, we move on the proof of Theorem \ref{the1.1}.

{\bf Proof of Theorem \ref{the1.1}}\quad  Let $\hat{\phi}\in \mathcal{C}^\infty_0(\mathbb{R})$ be an even, real-valued and non-negative funtion on $\R$ and satisfy
\begin{numcases}{\hat{\phi}(x)=}
1, &if $|x|\leq \frac{1}{4}$,\nonumber\\
0, &if $|x|\geq \frac{1}{2}$.\nonumber
\end{numcases}
Define the high frequency function $f_n$ and the low frequency function $g_n$ by
$$f_n=2^{-ns}\phi(x)\sin \bi(\frac{17}{12}2^nx\bi), \qquad g_n=2^{-\frac n2}\phi(x), \quad n\gg1.$$
It has been showed in \cite{Li 2020} that $\|f_n\|_{B_{p, r}^\sigma}\lesssim 2^{n(\sigma-s)}$.

Let
$$v^{1, n}(0, x)=(1-\partial_x)(f_n+g_n), \;v^{2, n}(0, x)=(1-\partial_x)f_n.$$
Consider Eq. (\ref{eq3.2}) with initial data $v^{1, n}(0, x)$ and $v^{2, n}(0, x)$, respectlively. Obviously, we have
\bbal
||v^{1, n}(0, x)-v^{2, n}(0, x)||_{B^{s-1}_{p,r}}=||(1-\partial_x)g_n||_{B^{s-1}_{p,r}}\leq C2^{-\frac{n}{2}},
\end{align*}
which means that
\bbal
\lim_{n\to\infty}||v^{1, n}(0, x)-v^{2, n}(0, x)||_{B^{s-1}_{p,r}}=0.
\end{align*}
It is easy to show that
\begin{eqnarray*}
||v^{1, n}(0, x)||_{B^{s-2}_{p,r}}&\lesssim& ||f_n+g_n||_{B^{s-1}_{p,r}}\lesssim ||f_n+g_n||_{B^{s-\frac{1}{2}}_{p,r}} \lesssim 2^{-\frac{n}{2}},\\
||v^{1, n}(0, x)||_{B^{s+\sigma}_{p,r}}&\lesssim& ||f_n+g_n||_{B^{s+\sigma+1}_{p,r}} \lesssim 2^{n(\sigma+1)} \qquad \mathrm{for} \qquad \sigma\geq -\frac{3}{2},\\
||v^{2, n}(0, x)||_{B^{s+t}_{p,r}}&\lesssim& ||f_n||_{B^{s+t+1}_{p,r}} \lesssim 2^{n(t+1)} \qquad \qquad\;\; \mathrm{for} \qquad t\in \mathbb{R},
\end{eqnarray*}
which imply
\bbal
&\big(||v^{1, n}(0, x)||^3_{B^{s-1}_{p,r}}+||v^{1, n}(0, x)||^2_{B^{s-2}_{p,r}}||v^{1, n}(0, x)||_{B^{s}_{p,r}}
+||v^{1, n}(0, x)||^{4}_{B^{s-2}_{p,r}}||v^{1, n}(0, x)||_{B^{s+1}_{p,r}}\big)\lesssim  1,
\\&\big(||v^{2, n}(0, x)||^3_{B^{s-1}_{p,r}}+||v^{2, n}(0, x)||^2_{B^{s-2}_{p,r}}||v^{2, n}(0, x)||_{B^{s}_{p,r}}
+||v^{2, n}(0, x)||^{4}_{B^{s-2}_{p,r}}||v^{2, n}(0, x)||_{B^{s+1}_{p,r}}\big)\lesssim  1.
\end{align*}
Furthermore, since $v^{1, n}(0, x)$ and $v^{2, n}(0, x)$ are both bounded in $B_{p, r}^{s-1}$, according to Proposition \ref{pro3}, we deduce that
\begin{eqnarray}\label{eq4.3}
&~&||\mathbf{S}_{t}(v^{1, n}(0, x))-\mathbf{S}_{t}(v^{2, n}(0, x))||_{B^{s-1}_{p,r}}\nonumber\\
&\geq&t\big|\big|([v^{1, n}(0, x)]^2-2u^{1, n}(0, x)v^{1, n}(0, x))\partial_xv^{1, n}(0, x)\nonumber\\
&~&-([v^{2, n}(0, x)]^2-2u^{2, n}(0, x)v^{2, n}(0, x))\partial_xv^{2, n}(0, x)\nonumber\\
&~&+\frac{1}{3}([v^{1, n}(0, x)]^3-[v^{2, n}(0, x)]^3)
+[\frac{1}{3}[u^{1, n}(0, x)]^3+\Phi_1(v^{1, n}(0, x))+
\Phi_2(v^{1, n}(0, x))]\nonumber\\
&~&-[\frac{1}{3}[u^{2, n}(0, x)]^3+\Phi_1(v^{2, n}(0, x))+
\Phi_2(v^{2, n}(0, x))]\big|\big|_{B^{s-1}_{p,r}}-\big|\big|(1-\partial_x)g_n\big|\big|_{B_{p, r}^{s-1}}-Ct^2\nonumber\\
&\geq&~ t\big|\big|([v^{1, n}(0, x)]^2-2u^{1, n}(0, x)v^{1, n}(0, x))\partial_xv^{1, n}(0, x)\nonumber\\
&~&-([v^{2, n}(0, x)]^2-2u^{2, n}(0, x)v^{2, n}(0, x))\partial_xv^{2, n}(0, x)\big|\big|_{B^{s-1}_{p, r}}-C2^{-\frac n2}-Ct^{2}.
\end{eqnarray}
For the sake of simplicity and convenience, in the following we denote
$$v^{i, n}(0, x)\triangleq v_i, \;\; u^{i, n}(0, x)\triangleq u_i, \;\; i=1, 2.$$
The coefficient of the first order term of $t$ in the last inequality in (\ref{eq4.3}) is simplified as
\begin{eqnarray*}
&\;&(v_1^2-2u_1v_1)\partial_xv_1-(v_2^2-2u_2v_2)\partial_xv_2\\
&=&(v_1^2\partial_xv_1-v_2^2\partial_xv_2)-2(u_1v_1\partial_xv_1-v_2\partial_xv_2)\\
&=&I_1-I_2.
\end{eqnarray*}
Bring in the concrete form of $v_1$ and $v_2$ where necessary, we have
\begin{eqnarray*}
I_1&=&(v_1^2-v_2^2)\partial_xv_1+v_2^2(\partial_xv_1-\partial_xv_2)\\
&=&(v_1^2-v_2^2)\partial_x(1-\partial_x)(f_n+g_n)+v_2^2(\partial_xv_1-\partial_xv_2)\\
&=&(v_1^2-v_2^2)\partial_x(f_n+g_n)-(v_1^2-v_2^2)\partial_x^2g_n-(v_1^2-v_2^2)\partial_x^2f_n+v_2^2\partial_x(v_1-v_2),\\
I_2&=&2u_1(v_1-v_2)\partial_xv_1+2v_2(u_1\partial_xv_1-u_2\partial_xv_2)\\
&=&2u_1(v_1-v_2)\partial_x(1-\partial_x)(f_n+g_n)+2v_2u_1\partial_x(v_1-v_2)+2v_2(u_1-u_2)\partial_xv_2\\
&=&2u_1(v_1-v_2)\partial_x(f_n+g_n)-2u_1(v_1-v_2)\partial_x^2g_n-2u_1(v_1-v_2)\partial_x^2f_n\\
&\;&+2v_2u_1\partial_x(v_1-v_2)+2v_2(u_1-u_2)\partial_xv_2.
\end{eqnarray*}
Using Lemma \ref{lem2.1}, after simple calculation, we obtain
\begin{eqnarray*}
\|(v_1^2-v_2^2)\partial_x(f_n+g_n)\|_{B_{p, r}^{s-1}}&\lesssim& \|v_1^2-v_2^2\|_{L^\infty}\|f_n+g_n\|_{B_{p, r}^s}+\|v_1^2-v_2^2\|_{B_{p, r}^{s-1}}\|\partial_x(f_n+g_n)\|_{L^\infty}\\
&\lesssim& 2^{-\frac{n}{2}},\\
\|(v_1^2-v_2^2)\partial_x^2g_n\|_{B_{p, r}^{s-1}}&\lesssim&\|v_1^2-v_2^2\|_{B_{p, r}^{s-1}}
\|\partial_x^2g_n\|_{B_{p, r}^{s-1}}\lesssim2^{-\frac{n}{2}},\\
\|v_2^2\partial_x(v_1-v_2)\|_{B_{p, r}^{s-1}}&\lesssim&\|v_2^2\|_{B_{p, r}^{s-1}}
\|\partial_x(v_1-v_2)\|_{B_{p, r}^{s-1}}\lesssim2^{-\frac{n}{2}},\\
\|2u_1(v_1-v_2)\partial_x(f_n+g_n)\|_{B_{p, r}^{s-1}}&\lesssim&\|u_1\|_{B_{p, r}^{s-1}}
\|(v_1-v_2)\|_{B_{p, r}^{s-1}}\|\partial_x(f_n+g_n)\|_{B_{p, r}^{s-1}}\lesssim2^{-\frac{n}{2}},\\
\|2u_1(v_1-v_2)\partial_x^2g_n\|_{B_{p, r}^{s-1}}&\lesssim&\|u_1\|_{B_{p, r}^{s-1}}
\|(v_1-v_2)\|_{B_{p, r}^{s-1}}\|\partial_x^2g_n\|_{B_{p, r}^{s-1}}\lesssim2^{-\frac{n}{2}},\\
\|2v_2u_1\partial_x(v_1-v_2)\|_{B_{p, r}^{s-1}}&\lesssim&\|v_2\|_{B_{p, r}^{s-1}}
\|u_1\|_{B_{p, r}^{s-1}}\|\partial_x(v_1-v_2)\|_{B_{p, r}^{s-1}}\lesssim2^{-\frac{n}{2}},\\
\|2v_2(u_1-u_2)\partial_xv_2\|_{B_{p, r}^{s-1}}&\lesssim&\|v_2\|_{B_{p, r}^{s-1}}
\|u_1-u_2\|_{B_{p, r}^{s-1}}\|\partial_xv_2\|_{B_{p, r}^{s-1}}\lesssim2^{-\frac{n}{2}}.
\end{eqnarray*}
While
\begin{eqnarray*}
[2u_1(v_1-v_2)-(v_1^2-v_2^2)]\partial_x^2f_n
&=&(1-\partial_x)g_n[g_n+\partial_x(2f_n+g_n)]\partial_x^2f_n\\
&=&(1-\partial_x)g_n(1+\partial_x)g_n\partial_x^2f_n+2(1-\partial_x)g_n\partial_xf_n\partial_x^2f_n,
\end{eqnarray*}
using product law (2) in Lemma \ref{lem2.1}, we have that
\begin{eqnarray*}
\|(1-\partial_x)g_n\partial_xf_n\partial_x^2f_n\|_{B_{p, r}^{s-1}}&\lesssim & \|(1-\partial_x)g_n\partial_xf_n\|_{B_{p, r}^{s-1}}\|\partial_x^2f_n\|_{L^\infty}+\|(1-\partial_x)g_n\partial_xf_n\|_{L^\infty}\|\partial_x^2f_n\|_{B_{p, r}^{s-1}}\\
&\lesssim & 2^{-\frac{n}{2}}\cdot2^0\cdot2^{n(2-s)}+2^{-\frac{n}{2}}\cdot2^{n(1-s)}\cdot2^n\lesssim 2^{-n(s-\frac{3}{2})}.
\end{eqnarray*}
Taking the above estimates into (\ref{eq4.3}), we find that
\begin{eqnarray}\label{eq4.4}
||\mathbf{S}_{t}(v^{1, n}(0, x))-\mathbf{S}_{t}(v^{2, n}(0, x))||_{B^{s-1}_{p,r}}
\geq t \big|\big|(1-\partial_x)g_n(1+\partial_x)g_n\partial_x^2f_n\big|\big|_{B^{s-1}_{p, r}}-C2^{-\frac n2}-Ct^{2}.
\end{eqnarray}
For the term $(1-\partial_x)g_n(1+\partial_x)g_n\partial_x^2f_n$, it can be verified that $\Delta_j\big((1-\partial_x)g_n(1+\partial_x)g_n\partial_x^2f_n\big)=0,$  $j\neq n$ and $\Delta_n\big((1-\partial_x)g_n(1+\partial_x)g_n\partial_x^2f_n\big)=(1-\partial_x)g_n(1+\partial_x)g_n\partial_x^2f_n$
for $n\geq 5.$ Direct calculation shows that
\begin{eqnarray*}
&\;&\big|\big|(1-\partial_x)g_n(1+\partial_x)g_n\partial_x^2f_n\big|\big|_{B^{s-1}_{p, r}}=2^{n(s-1)}\big|\big|(1-\partial_x)g_n(1+\partial_x)g_n\partial_x^2f_n\|_{L^p}\\
&=&\|2^{-2n}(1-\partial_x)\phi(1+\partial_x)\phi\partial_x^2\phi \sin(\frac{17}{12}2^nx)+\frac{17}{12}2^{-n}(1-\partial_x)\phi(1+\partial_x)\phi\phi \cos(\frac{17}{12}2^nx)\\
&\;&+\frac{17}{12}2^{-n}(1-\partial_x)\phi(1+\partial_x)\phi\partial_x\phi \cos(\frac{17}{12}2^nx)
+\frac{17}{12}(1-\partial_x)\phi(1+\partial_x)\phi\partial_x\phi \sin(\frac{17}{12}2^nx)
\big|\big|_{L^p}\\
&\geq&\frac{17}{12}\big|\big|\psi(x)\sin(\frac{17}{12}2^nx)\big|\big|_{L^p}-2^{-n}\rightarrow \frac{17}{12}\big(\frac{\int_0^{2\pi}|\sin x|^pdx}{2\pi}\big)^{\frac{1}{p}}\big|\big|\psi(x)\big|\big|_{L^p} ,
\end{eqnarray*}
by the Riemann Theorem, where $(1-\partial_x)\phi(1+\partial_x)\phi\partial_x\phi \triangleq \psi(x),$
which together with (\ref{eq4.4}) yield
\bbal
\liminf_{n\rightarrow \infty}||\mathbf{S}_{t}(v^{1, n}(0, x))-\mathbf{S}_{t}(v^{2, n}(0, x))||_{B^{s-1}_{p,r}}\gtrsim t\quad\text{for} \ t \ \text{small enough}.
\end{align*}
That is to say, the solution map $v_0\rightarrow S_t(v_0)$ of the initial value problem (\ref{eq3.2}) depends not uniformly continuous on initial data in $B_{p, r}^{s-1}$.

Since $u=(1-\partial_x)^{-1}v$ and $(1-\partial_x)^{-1}$ is an isomorphic mapping from $B_{p, r}^{s-1}(\mathbb{R})$ into $B_{p, r}^s(\mathbb{R})$, hence the non-uniform continuous dependence of $v$ in $B_{p, r}^{s-1}$ is consistent with that of $u$ in $B_{p, r}^s.$

Thus, this completes the proof of Theorem \ref{the1.1}.

\section*{Acknowledgments}
The author is very grateful to Dr. Jinlu Li for some
useful suggestions. This work is partially supported by the National Natural Science Foundation of China (Grant No.11801090).
%\vspace*{1em}

\end{document}